\documentclass[12pt,a4paper]{amsart}

\usepackage{amssymb}

\newtheorem*{Def}{Definition}
\expandafter\def\expandafter\Def\expandafter{\Def\upshape}
\newtheorem*{Question}{Question}

\newtheorem{Result}{Result}[section]
\newtheorem{Thm}  [Result]{Theorem}
\newtheorem{Prop} [Result]{Proposition}
\newtheorem{Lemma}[Result]{Lemma}
\newtheorem{Cor}  [Result]{Corollary}
\newtheorem{Remark}[Result]{Remark}
\expandafter\def\expandafter\Remark\expandafter{\Remark\upshape}

\def\qedbox{\hbox{\vrule\vbox{\hrule width 6pt\vskip6pt\hrule}\vrule}}
\def\qed{\ifvmode\leavevmode\fi
  \unskip\nobreak\hfill\penalty50 \quad \null\nobreak\hfill
  \qedbox{\parfillskip0pt \finalhyphendemerits0 \par}}
\newenvironment{Proof}{\ifdim\lastskip=0pt \medskip\fi
  \noindent{\it Proof.~}\ignorespaces}{\qed\medskip}

\def\eqalign#1{\null\,\vcenter{\openup\jot\mathsurround0pt\relax
  \ialign{\strut\hfil$\displaystyle{##}$&$\displaystyle{{}##}$\hfil
      \crcr#1\crcr}}\,}

\def\N{{\mathbb{N}}}
\let\eps\varepsilon

\let\<\langle
\let\>\rangle

\let\maparrow\longrightarrow
\def\map#1#2#3{#1\colon#2\,{\maparrow}\,#3}

\def\vv{<\!\!<}

\catcode`\@=11
\let\tendstoarrow\longrightarrow
\def\tendst@style#1#2#3#4{{\setbox0\hbox{$#1\tendstoarrow$}%
  \dimen0 \ht0 \advance\dimen0 -\fontdimen22#22
  #1\mathop{\vcenter to 2\dimen0{\box0\vss}}\limits#3#4}}
\def\tendst@#1#2{\mathrel{\mathchoice
 {\tendst@style\displaystyle\textfont{#1}{#2}}
 {\tendst@style\textstyle\textfont{#1}{#2}}
 {\tendst@style\scriptstyle\scriptfont{#1}{#2}}
 {\tendst@style\scriptscriptstyle\scriptscriptfont{#1}{#2}}
}}
\newtoks\tendst@toks
\def\tendst@none{\tendst@{}{}}
\def\tendst@one{\tendst@{\the\tendst@toks}{}}
\def\tendst@pw#1#2{\tendst@{\the\tendst@toks}{_{#2}}}
\def\tendst@bw#1#2{\tendst@{\the\tendst@toks}{^{#2}}}
\def\tendst@sup{\relax \ifcat_\noexpand\next \expandafter\tendst@pw \else
 \expandafter\tendst@one \fi}
\def\tendst@sub{\relax \ifcat^\noexpand\next \expandafter\tendst@bw \else
 \expandafter\tendst@one \fi}
\def\tendst@setsup#1#2{\tendst@toks{^{#2}}\futurelet\next\tendst@sup}
\def\tendst@setsub#1#2{\tendst@toks{_{#2}}\futurelet\next\tendst@sub}
\def\tendst@s{\relax \ifcat^\noexpand\next \expandafter\tendst@setsup \else
 \ifcat_\noexpand\next \expandafter\expandafter\expandafter\tendst@setsub
 \else \tendst@none \fi\fi}
\def\tendsto{\futurelet\next\tendst@s}
\catcode`\@=12

\def\LIM{{\rm LIM}}

\def\Jphi{{J(\Phi)}}
\def\JJ{{\hat J}}
\def\JJphi{{\JJ(\Phi)}}

\let\psf\wp

\begin{document}

\title{Complemented subspaces of $J$-sums of Banach spaces}

\author{Manuel Gonz\'alez}
\author{Javier Pello}

\begin{abstract}
We study the complemented subspaces of the $J$-sums of Banach spaces $J(\Phi)$
and $\hat J(\Phi)$ introduced by Bellenot. As an application, we show that,
under some conditions, $J(\Phi)$ and $\hat J(\Phi)$ are subprojective,
i.e., every closed infinite-dimensional subspace of either of them
contains a complemented infinite-dimensional subspace.
\end{abstract}

\maketitle

\begin{section}{Introduction}

Extending a result of James \cite{james:60}, Lindenstrauss \cite{linden:71}
showed that for every separable Banach space~$X$ there exists a separable
Banach space Y with a monotone shrinking basis and a surjective operator~$Q$
from~$Y^*$ onto~$X$ such that $Y^{**} = i_Y(Y)\oplus Q^*(X^*)$, where
$\map{i_Y}{Y}{Y^{**}}$ is the canonical embedding (hence $Y^{**}\simeq
Y\oplus X^*$), and the same result was obtained for $X$ and~$Y$
weakly compactly generated in \cite{DFJP}.

Bellenot gave a general construction in \cite{Jsum} subsuming both results.
Given  a sequence of Banach spaces $(X_n, \|\cdot\|_n)_{n\in\N_0}$ with
$\dim X_0 = 0$ and a sequence of operators $\Phi = (\phi_n)_{n\in\N_0}$ with
$\map{\phi_n}{X_n}{X_{n+1}}$ and $\|\phi_n\| \leq 1$ for every $n\in\N_0$,
and denoting $\Pi = \prod_{n\in\N_0} X_n$, he defined a quantity $\|x\|_J$
for every $x = (x_n)_{n\in\N_0} \in \Pi$ that is similar to the norm in the
classical James space $J$, and considered the Banach spaces
$\JJphi = \{ x \in \Pi : \|x\|_J < \infty \}$ and
$\Jphi = \{ x \in \JJphi : \lim_n \|x_n\|_n = 0 \}$.


Here we study some families of projections in the spaces $\JJphi$\break
and~$\Jphi$ and prove that, if each $X_n$ is subprojective, then every closed
infinite-dimensional subspace of~$\Jphi$ contains an infinite-dimensional
subspace complemented in~$\JJphi$, so $\Jphi$ itself is subprojective, and
if furthermore the quotient $\JJphi/\Jphi$ is subprojective or the quotient
map $\JJphi \maparrow \JJphi/\Jphi$ is strictly singular, then also $\JJphi$
is subprojective.
We also give conditions for $\Jphi$ to be complemented or not 
in~$\JJphi$ and describe several examples illustrating the scope of these
results.


An operator $\map{T}{X}{Y}$ between two Banach spaces $X$ and~$Y$ is called
\emph{strictly singular} if there is no closed infinite-dimensional
subspace~$M$ of~$X$ such that the restriction $T|_M$ is an isomorphism.
A Banach space~$X$ is called \emph{subprojective} if every closed
infinite-dimensional subspace of~$X$ contains an infinite-dimensional
subspace complemented in~$X$; note that finite-dimensional spaces are
trivially subprojective. Subprojective spaces were introduced by
Whitley~\cite{whitley} to find conditions for the conjugate of an operator
to be strictly singular. 
There has been a surge in attention to these spaces after a recent
systematic study~\cite{oikhberg-spinu}, such as to obtain some positive
solutions to the perturbation classes problem for semi-Fredholm operators,
which has a negative solution in general~\cite{gonzalez} but there are some
positive answers when one of the spaces is subprojective 
\cite{gonzalez-et-al} \cite{G-P-Salas}.

\end{section}

\begin{section}{Definitions and basic facts}

Let $\N_0 = \N\cup\{0\}$ and
$\psf = \{\, S \subset \N_0 : \hbox{$S$ non-empty, finite} \,\}$.
In the sequel, $(X_n, \|\cdot\|_n)_{n\in\N_0}$ is a sequence of Banach spaces
with $\dim X_0 = 0$ and $\Phi = (\phi_n)_{n\in\N_0}$ is a sequence of
operators $\map{\phi_n}{X_n}{X_{n+1}}$ such that $\|\phi_n\| \leq 1$
for every $n\in\N_0$. If $n\leq m\in\N_0$, we will write $\map{\phi_n^m =
\phi_{m-1}\circ\cdots\circ\phi_{n+1}\circ\phi_n}{X_n}{X_m}$,
with the convention that $\phi_n^n = I_{X_n}$ is the identity on~$X_n$;
in particular, $\phi_n^{n+1} = \phi_n$ for every $n\in\N_0$ and
$\phi_n^{m+1} = \phi_m \circ \phi_n^m$ for every $n\leq m\in\N_0$.

Denote $\Pi = \prod_{n\in\N_0} X_n$ and,
given $x = (x_n)_{n\in\N_0} \in \Pi$ and
$S = \{\, p_0 < \cdots < p_k \,\} \in \psf$, define
$$\eqalign{
  \sigma(x, S) &=
    \biggl( \sum_{i=1}^k \bigl\| \phi_{p_{i-1}}^{p_i}(x_{p_{i-1}})
      - x_{p_i} \bigr\|_{p_i} ^2\biggr)^{1/2} \cr
  \rho(x, S) &=
    \biggl( \sum_{i=1}^k \bigl\| \phi_{p_{i-1}}^{p_i}(x_{p_{i-1}})
      - x_{p_i} \bigr\|_{p_i} ^2 + \|x_{p_k}\|_{p_k}^2 \biggr)^{1/2} \cr
}$$
and $\|x\|_J = (1/\sqrt2) \sup_{S\in\psf} \rho(x,S)$.
Note that $\sigma(\,\cdot\,,S)$ and $\rho(\,\cdot\,,S)$ are seminorms
in~$\Pi$ and that $\rho(x,S)^2 = \sigma(x,S)^2 + \|x_{p_k}\|_{p_k}^2$, so
$\sigma(x,S) \leq \rho(x,S)$. Also, $\sigma(x,S) = 0$ if $S$ has a single
element, $\sigma(x, S\cup T)^2 = \sigma(x,S)^2 + \sigma(x,T)^2$ if $\max S =
\min T$ and $\sigma(x, S\cup T)^2 \geq \sigma(x,S)^2 + \sigma(x,T)^2$
if $\max S \leq \min T$.

\begin{Def}
\cite{Jsum}
The $J$-sum $\Jphi$ is defined as the completion of the normed space of
the finitely non-zero sequences in~$\Pi$ with~$\|\cdot\|_J$.
\end{Def}

\begin{Prop}
\label{bellenot-J}
\cite[p.~98, Remarks 3 and~4]{Jsum}
Each $X_m$ can be identified isometrically with the subspace of sequences
$(x_n)_{n\in\N_0} \in \Jphi$ such that $x_n = 0$ for $n\ne m$.
With this identification, $(X_n)_{n\in\N}$ is a bimonotone decomposition
for~$\Jphi$, 
i.e., $\bigl\| \sum_{i=p+1}^{p+q} x_i \bigr\|_J \leq
\bigl\| \sum_{i=1}^{p+q+r} x_i \bigr\|_J$ for every
$(x_n)_{n\in\N_0} \in \Jphi$ and $p$, $q$, $r\in\N_0$.
\end{Prop}

In particular, $\|(x_n)_{n\in\N_0}\|_J =
\sup_{n\in\N} \bigl\| \sum_{i=1}^n x_i \bigr\|_J =
\lim_n \bigl\| \sum_{i=1}^n x_i \bigr\|_J$ for every $(x_n)_{n\in\N_0}
\in \Jphi$.

Define now
$$\JJphi = \bigl\{\, x \in \Pi : \|x\|_J < \infty \,\bigr\},$$
which can be identified with $(X_n)_{n\in\N}^\LIM$, the set of all sequences
$(x_n)_{n\in\N}\in\Pi$ such that $(\sum_{i=1}^n x_i)_{n\in\N}$ is bounded
\cite[p.~96]{Jsum} \cite[Proposition~1.b.2]{classical1}. This identification
makes for slightly simpler notation than the direct use of
$(X_n)_{n\in\N}^\LIM$.

The space of eventually constant sequences is
$$\eqalign{
  \Omega(\Phi) = \{\, &(x_n)_{n\in\N_0} \in \Pi : \cr
  &\quad
    \hbox{there is $n\in\N$ such that $\phi_m(x_m) = x_{m+1}$ for all $m\geq n$}
  \,\}. \cr
}$$
Clearly $\Omega(\Phi) \subset \JJphi$ and $\|(x_n)_{n\in\N_0}\|_\Omega =
\lim_n \|x_n\|_n$ defines a seminorm on~$\Omega(\Phi)$, as $\|\phi_n\|\leq1$
for every $n\in\N_0$. We denote the completion of the normed space
$(\Omega(\Phi)/\mathop{\rm ker}\|\cdot\|_\Omega,\|\cdot\|_\Omega)$
by~$\tilde\Omega(\Phi)$.

\begin{Thm}
\label{bellenot-JJ}
\cite[Theorem~1.1]{Jsum}
\begin{itemize}
\item[(i)] $\Omega(\Phi)$ is dense in~$\JJphi$.
\item[(ii)] The unique extension $\map\Theta\JJphi{\tilde\Omega(\Phi)}$ of
the natural map $\Omega(\Phi) \maparrow \tilde\Omega(\Phi)$ is a quotient map
with kernel~$\Jphi$.
\item[(iii)] If each $X_n$ is reflexive, then $\Jphi^{**} = \JJphi$, hence
$\Jphi^{**}/\Jphi$ is isometric to~$\tilde\Omega(\Phi)$.
\end{itemize}
\end{Thm}

In particular, $\Jphi$ can be identified with
$$\bigl\{\, (x_n)_{n\in\N_0} \in \JJphi : \|x_n\|_n \tendsto_n 0 \,\bigr\}.$$

\begin{Remark}
\label{not-2}
It is worth noting that the choice of the $\ell_2$-norm in the definition of
$\sigma$ and~$\rho$ above is not essential to the construction of $\JJphi$;
all of the results below hold equally well for any $\ell_q$-norm
with $1 < q < \infty$.
\end{Remark}

\end{section}

\begin{section}{Estimates for disjoint sequences}

We will adopt the following definition: If $A$, $B$ are non-empty subsets
of~$\N_0$, we will write $A \vv B$ if there exists $n\in\N$ such that
$a < n < b$ for every $a\in A$ and $b\in B$; equivalently,
$\min B - \max A \geq 2$.

\begin{Def}
Let $I\subseteq\N_0$ be an interval. We define $\map{P_I}{\JJphi}{\JJphi}$ as
$$P_I \bigl( (x_n)_{n\in\N_0} \bigr) = ( x_n \chi_I(n) )_{n\in\N_0}.$$
\end{Def}

Each $P_I$ is well defined and a projection with $\|P_I\| \leq 1$ due to the
following result with $m = 0$.

\begin{Prop}
\label{J-interval}
Let $I\subseteq\N$ be an interval, let $x = (x_n)_{n\in\N_0} \in \Pi$ such
that $x_n = 0$ for every $n\notin I$, let $m < \min I$ and let $S\in\psf$.
Then $\rho(x,S) \leq \rho(x,\{m\}\cup(S\cap I))$.
\end{Prop}

\begin{Proof}
Let $\tilde S = S \cap [\min I,\infty) \subseteq S$; then $\rho(x,S) \leq
\rho(x,\{m\}\cup\tilde S)$. Indeed, if $\tilde S \subsetneq S$, then the
elements in $S \setminus \tilde S$ do not contribute to $\rho(x,S)$ except
for the possible effect of their presence on the first element of~$\tilde S$,
and then $\{m\}$ has the same effect; and, if $S = \tilde S$, then the
addition of $\{m\}$ only prepends a (non-negative) term to $\rho(x,S)$.

Now, the elements in $\tilde S$ beyond the maximum of~$I$, if any,
do not contribute to $\rho(x,S)$ either, except for the fact that the last
term in~$\rho(x,S)$ may become smaller (as every $\|\phi_n\| \leq 1$).
In any case, $\rho(x,S) \leq \rho(x,\{m\}\cup\tilde S) \leq
\rho(x,\{m\}\cup(S\cap I))$.
\end{Proof}

As a consequence, if $I\subset\N_0$ is a finite interval and $x\in R(P_I)$,
then $\sup_{S\in\psf} \rho(x,S)$ is attained for some $S\subseteq[0,\max I]$.

If $n\in\N_0$, we define the projection $P_n = P_{\{1,\ldots,n\}}$
(as opposed to $P_{\{n\}}$, which gives the identification of~$X_n$ described
in Proposition~\ref{bellenot-J}); in particular, $P_0 = P_\emptyset = 0$.

Sequences in~$\JJphi$ with disjoint supports admit an upper $2$-estimate,
and also a lower $2$-estimate if the supports are furthermore not adjacent.

\begin{Prop}
\label{J-upper}
\cite[Theorem 1.1(III)]{Jsum}
Let $I_1$, $I_2$, \dots, $I_m \subseteq \N_0$ be disjoint intervals
and let $x_j \in R(P_{I_j})$ for every $j\in\{1,2,\ldots,m\}$.
Then $\| \sum_{j=1}^m x_j \|_J^2 \leq 3 \sum_{j=1}^m \|x_j\|_J^2$.
\end{Prop}

\begin{Prop}
\label{J-lower}
Let $I_1 \vv I_2 \vv \cdots \vv I_m \subseteq \N_0$ be non-empty\break
intervals and let $x_j \in R(P_{I_j})$ for every $j\in\{1,2,\ldots,m\}$.
Then\break $2 \| \sum_{j=1}^m x_j \|_J^2 \geq \sum_{j=1}^m \|x_j\|_J^2$.
\end{Prop}

\begin{Proof}
Define $y = (y_n)_{n\in\N_0} = \sum_{j=1}^m x_j \in \JJphi$, let $\eps > 0$
and $n_1 = 0$ and pick $n_j \in (\max I_{j-1}, \min I_j)$ for every
$j\in\{2,\ldots,m\}$, which is possible because $I_{j-1} \vv I_j$; note
that $y_{n_j} = 0$ for all $j\in\{1,2,\ldots,m\}$.

For every $j\in\{1,\ldots,m\}$, there exists $S_j\in\psf$ such that
$\rho(x_j,S_j)^2 > 2 \|x_j\|_J^2 - \eps/m$ and,
by Proposition~\ref{J-interval}, we may assume $n_j \in S_j \subseteq \{n_j\}
\cup I_j$. Let then $m_j = \max S_j$, so 
$$\rho(y,S_j)^2 = \sigma(y,S_j)^2 + \|y_{m_j}\|_{m_j}^2 =
  \sigma(y,S_j)^2 + \sigma(y,\{n_j,m_j\})^2$$
and there exists some $\tilde S_j$ (either $S_j$ itself or $\{n_j,m_j\}$)
such that\break $\sigma(y,\tilde S_j)^2 \geq {1\over2} \rho(y,S_j)^2$, with
$\max\tilde S_{j-1} = m_{j-1} < n_j = \min\tilde S_j$ if $j > 1$.
Let $\tilde S = \bigcup_{j=1}^m \tilde S_j$; then
$$\eqalign{
  2 \|y\|_J^2 & \geq \rho(y,\tilde S)^2 \geq \sigma(y,\tilde S)^2
      \geq \sum_{j=1}^m \sigma(y,\tilde S_j)^2
      \geq {1\over2} \sum_{j=1}^m \rho(y,S_j)^2 \cr
    & = {1\over2} \sum_{j=1}^m \rho(x_j,S_j)^2
      > {1\over2} \sum_{j=1}^m \bigl( 2 \|x_j\|_J^2 - \eps/m \bigr)
      > \sum_{j=1}^m \|x_j\|_J^2 - \eps. \cr
}$$
As this is true for all $\eps > 0$, it follows that $2\|y\|_J^2 \geq
\sum_{j=1}^m \|x_j\|_J^2$ indeed.
\end{Proof}

\end{section}

\break

\begin{section}{Stepping projections}

We will write $\Lambda$ for the set of all sequences $(\alpha_n)_{n\in\N_0}$
in~$\N_0$ such that $\alpha_n \leq \alpha_{n+1}$ and $\alpha_n \leq n$ for
every $n\in\N_0$ and, given $\alpha = (\alpha_n)_{n\in\N_0} \in \Lambda$,
we define a linear map $\map{Q_\alpha}{\Pi}{\Pi}$ as
$$Q_\alpha \bigl( (x_n)_{n\in\N_0} \bigr) =
  \bigl( \phi_{\alpha_n}^n (x_{\alpha_n}) \bigr) _{n\in\N_0}.$$

\begin{Lemma}
\label{J-Q-S}
Let $\alpha = (\alpha_n)_{n\in\N_0} \in \Lambda$, let $x \in \JJphi$
and let $S\in\psf$. Then
$$\rho(Q_\alpha(x),S) \leq \rho(x, \{\, \alpha_p : p\in S \,\}).$$
\end{Lemma}

\begin{Proof}
Write $x = (x_n)_{n\in\N_0} \in \JJphi$ and let
$y_n = \phi_{\alpha_n}^n (x_{\alpha_n})$ for every $n\in\N_0$,
so that $Q_\alpha(x) = (y_n)_{n\in\N_0}$.

Let $p < q$ be two consecutive elements of~$S$. Then
$$\eqalign{
  \| \phi_p^q (y_p) - y_q \|_q
    &= \| \phi_p^q (\phi_{\alpha_p}^p (x_{\alpha_p}))
        - \phi_{\alpha_q}^q (x_{\alpha_q}) \|_q \cr
    &= \| \phi_{\alpha_q}^q (\phi_{\alpha_p}^{\alpha_q} (x_{\alpha_p}))
        - \phi_{\alpha_q}^q (x_{\alpha_q}) \|_q \cr
    &= \| \phi_{\alpha_q}^q (\phi_{\alpha_p}^{\alpha_q} (x_{\alpha_p})
        - x_{\alpha_q}) \|_q \cr
    &\leq \| \phi_{\alpha_p}^{\alpha_q} (x_{\alpha_p})
        - x_{\alpha_q} \|_{\alpha_q}; \cr
}$$
in particular, $\| \phi_p^q (y_p) - y_q \|_q = 0$ if $\alpha_p = \alpha_q$.
Additionally, for the final term of $\rho(Q_\alpha(x),S)$, if $p = \max S$,
we also have
$$\|y_p\|_p = \| \phi_{\alpha_p}^p (x_{\alpha_p}) \|_p \leq
  \|x_{\alpha_p}\|_{\alpha_p},$$
so indeed $\rho(Q_\alpha(x),S) \leq \rho(x, \{\, \alpha_p : p\in S \,\})$.
\end{Proof}

\begin{Prop}
Let $\alpha\in\Lambda$. Then $\map{Q_\alpha|_\JJphi}{\JJphi}{\JJphi}$
is a bounded linear operator and $\|Q_\alpha\|\leq1$.
\end{Prop}

\begin{Proof}
$Q_\alpha$ is clearly linear and $\|Q_\alpha(x)\|_J \leq \|x\|_J$
for every $x \in \JJphi$ by Lemma~\ref{J-Q-S}.
\end{Proof}

\begin{Def}
\upshape
Let $A\subseteq\N$. We define a projection $\map{Q_A}{\JJphi}{\JJphi}$ as
$Q_A = Q_{\alpha_A}$ where $\alpha_A = (\alpha_n)_{n\in\N_0}$ is given by
$\alpha_n = \max ( (\{0\} \cup A) \cap [0,n] )$.
If $n\in\N_0$, we define the projection $Q_n = Q_{\{1,\ldots,n\}}$
(as opposed to $Q_{\{n\}}$); in particular, $Q_0 = Q_\emptyset = 0$.
\end{Def}

If $A = \{ a_1 < a_2 < \ldots \} \subseteq \N$, and $A_0 = \{0\} \cup A$,
then $Q_A((x_n)_{n\in\N_0}) = (y_m)_{m\in\N_0}$ is given by $y_m =
\phi_a^m(x_a)$, where $a = \max A_0\cap[0,m]$, for every $m\in\N$,
so $Q_A(x)$ depends only on those components $x_n$ for which $n\in A$, and
$N(Q_A) = \{\, (x_n)_{n\in\N_0} : \hbox{$x_n=0$ for every $n\in A$} \,\}$.
If $(I_j)_{j\in\N}$ is a sequence of non-empty intervals of~$\N_0$ such that
$I_j \vv I_{j+1}$ for every $j\in\N$ and
$A = \N\setminus\bigcup_{j\in\N} I_j$, then $N(Q_A)$
is isomorphic to $\ell_2((R(P_{I_j}))_{j\in\N})$ in the natural way
due to Propositions \ref{J-upper} and~\ref{J-lower}.

If $x = (x_m)_{m\in\N_0} \in \JJphi$ and $n\in\N$, then
$$Q_n(x) = 
  (x_1, x_2, \ldots, x_{n-1}, x_n,
    \phi_n^{n+1}(x_n), \phi_n^{n+2}(x_n), \ldots),$$
so $Q_n$ is a projection (clearly $Q_n^{\;2} = Q_n$) and the set of
eventually constant sequences can be written as $\Omega(\Phi) =
\bigcup_{n\in\N} R(Q_n)$, which is dense in~$\JJphi$,
but this can be stated in a better form.

\begin{Prop}
\label{J-Qlim}
Let $x \in \JJphi$. Then $Q_n(x) \tendsto_n x$.
\end{Prop}

\begin{Proof}
The result is immediate for the set of eventually constant sequences
$\Omega(\Phi) = \bigcup_{n\in\N} R(Q_n)$, which is dense in~$\JJphi$
by Theorem~\ref{bellenot-JJ}.
\end{Proof}

Given $n\in\N$, both $P_n$ and~$Q_n$ are projections that depend only on
the first $n$ components of~$x$, so they have the same kernel, and their
ranges are different but isometric.

\begin{Prop}
\label{J-PQ-QP}
Let $m\in\N$. Then $Q_mP_m = Q_m$ and $P_mQ_m = P_m$ and, for every
$x\in\JJphi$, $\|P_m(x)\|_J = \|Q_m(x)\|_J$.
\end{Prop}

\begin{Proof}
The identities $Q_mP_m = Q_m$ and $P_mQ_m = P_m$ follow from the fact that
both $P_m(x)$ and $Q_m(x)$ depend only on the first $m$ components of~$x\in
\JJphi$ and leave them untouched, and then
$$\|Q_m(x)\|_J = \|Q_mP_m(x)\|_J \leq \|P_m(x)\|_J
               = \|P_mQ_m(x)\|_J \leq \|Q_m(x)\|_J$$
from $\|P_m\| \leq 1$ and $\|Q_m\| \leq 1$.
\end{Proof}

\end{section}

\begin{section}{Subspaces of $\Jphi$ and~$\JJphi$}

As already mentioned in Theorem~\ref{bellenot-JJ}, $\Jphi^{**}$ can be
identified with $\JJphi$ when every $X_n$ is reflexive. However, even if
this is not the case, $\Jphi^{**}$ can be identified isometrically with
$\JJ(\Phi^{**})$, where $\Phi^{**} = (\phi_n^{**})_{n\in\N_0}$,
\cite[p.~97]{Jsum}, and then $\JJphi$ can be seen to embed isometrically
into $\JJ(\Phi^{**}) \equiv \Jphi^{**}$ by way of the natural inclusion
of each $X_n$ into~$X_n^{**}$, where each $x\in\JJphi$ is identified with
the weak$^*$ limit of $(P_n(x))_{n\in\N}$ in~$\Jphi^{**}$
(if every $X_n$ is reflexive, then simply $\Phi^{**} = \Phi$).

\begin{Prop}
\label{subproj-helper}
Let $Y$ be a subprojective Banach space, let\break $\map{T}{X}{Y}$ be an operator
and let $M$ be a closed infinite-dimensional subspace of~$X$ such that $T|_M$
is not strictly singular. Then $M$ contains an infinite-dimensional subspace
complemented in~$X$.
\end{Prop}

\begin{Proof}
$T|_M$ is not strictly singular, so there exists some
infinite-\break dimensional
subspace~$N$ of~$M$ such that $T|_N$ is an isomorphism. As $Y$ is
subprojective, we can further refine~$N$ to assume that $T(N)$ is
complemented in~$Y$, and then $Y = T(N) \oplus Z$ implies
$X = N \oplus T^{-1}(Z)$, so $N$ is complemented in~$X$.
\end{Proof}

\break

\begin{Thm}
\label{J-main}
(i) If each $X_n$ is subprojective, then every closed infinite-dimensional
subspace of~$\Jphi$ contains an infinite-dimensional subspace complemented
in~$\JJphi$;
\hfil\break\indent
(ii) if each $X_n$ is hereditarily~$\ell_2$, then $\Jphi$
is hereditarily~$\ell_2$;
\hfil\break\indent
(iii) if each $X_n$ is subprojective and hereditarily~$\ell_2$, then every
closed infinite-dimensional subspace of~$\Jphi$ contains a copy of~$\ell_2$
complemented in~$\JJphi$.
\end{Thm}

Note that case~(ii) is already hinted at in \cite{Jsum} and proved
in~\cite[Lemma~3]{ostrovskii} when each $X_n$ is finite-dimensional.

\begin{Proof}
Let $M$ be a closed infinite-dimensional subspace of~$\Jphi$.
If $Q_n|_M$ is not strictly singular for some $n\in\N$, then $R(Q_n)$ is
isometric to~$R(P_n)$ by Proposition~\ref{J-PQ-QP}, in turn isomorphic to
$\bigoplus_{i=1}^n X_i$, so, respectively for each case,

(i) $M$ contains an infinite-dimensional subspace complemented\break
in~$\JJphi$ by Proposition~\ref{subproj-helper}, since $R(Q_n)$ is
isomorphic to $\bigoplus_{i=1}^n X_i$, which is subprojective
\cite[Proposition~2.2]{oikhberg-spinu}; or

(ii) $M$ contains an infinite-dimensional subspace~$N$ such that $Q_n|_N$
is an isomorphism, where $R(Q_n)$ is isomorphic to $\bigoplus_{i=1}^n X_i$,
so $N$ contains a copy of~$\ell_2$ as containing a copy of~$\ell_2$ is a
three-space property~\cite{3space}; or

(iii) both of the above apply, so $M$ contains a copy of~$\ell_2$ and then
said copy contains a further copy of~$\ell_2$ complemented in~$\JJphi$.

Otherwise, assume that $Q_n|_M$ is strictly singular for every $n\in\N$,
in any of the cases.
Then, for every $n\in\N$ and $\eps > 0$, there exists $x\in M$ such that
$\|x\| = 1$ and $\|Q_n(x)\| < \eps$, and then there is $m > n + 1$ such
that $\|P_{m-1}(x) - x\| < \eps$.
By induction, starting with an arbitrary $n_1\in\N$, there exist a sequence
$(n_k)_{k\in\N}$ of elements in~$\N$ and a sequence $(x_k)_{k\in\N}$ of
norm-one elements in~$M$ such that $\|Q_{n_k}(x_k)\| < 2^{-k}/6$,
$n_{k+1} > n_k + 1$ and $\|P_{n_{k+1}-1}(x_k) - x_k\| < 2^{-k}/6$
for every $k\in\N$.

For each $k\in\N$, define $T_k = P_{n_{k+1}-1} (I - Q_{n_k})$ and note that,
given $y = (y_n)_{n\in\N_0} \in \JJphi$ and $m\in\N$, then 
$$T_k(y)_m = \left\{\,\vcenter{\normalbaselines\openup\jot
  \halign{$#$\hfil&\quad#\hfil\cr
    y_m - \phi_{n_k}^m(y_{n_k}),& if $n_k < m < n_{k+1}$ \cr
    0,& otherwise \cr
}}\right.$$
so $T_k(y)$ depends only on components $[n_k,n_{k+1})$ of~$y$ and its value
lies in the range of $P_{(n_k,n_{k+1})}$. As a consequence, it is easy to
check that each $T_k$ is a projection with $\|T_k\| \leq 2$
and that $T_iT_j = 0$ if $i \ne j$.
If we further define $A = \{\, n_k : k\in\N \,\}$ and $I_k = (n_k,n_{k+1})$
for each $k\in\N$, it holds that
$$T_k = P_{n_{k+1}-1} (I - Q_{n_k}) = P_{I_k} (I - Q_{n_k})
  = P_{I_k} (I - Q_A)$$
for every $k\in\N$.

Let $z_k = T_k(x_k) = P_{n_{k+1}-1} (I - Q_{n_k}) (x_k) \in R(P_{I_k})$ for
every $k\in\N$; then $\|z_k - x_k\| \leq \| P_{n_{k+1}-1}(x_k) - x_k \| +
\| P_{n_{k+1}-1} Q_{n_k} (x_k) \| < 2^{-k}/6 + 2^{-k}/6 = 2^{-k}/3 \leq 1/6$,
so $5/6 < \|z_k\| < 7/6$ for every $k\in\N$. If we take $x^*_k \in \Jphi^*$
such that $\|x^*_k\| < 6/5$ and $\<x^*_k, z_k\> = 1$ for each $k\in\N$,
then, for every $x\in\JJphi$ and $k\in\N$, Proposition~\ref{J-lower} yields
$$\eqalign{
  \sum_{i=1}^k \bigl| \< x^*_i, T_i(x) \> \bigr| ^2
    & \leq \sum_{i=1}^k 2 \|T_i(x)\|_J ^2
      \leq  4 \biggl\| \sum_{i=1}^k T_i(x) \biggr\|_J ^2 \cr
    & =    4 \biggl\| \sum_{i=1}^k P_{I_i} (I - Q_A) (x) \biggr\|_J ^2 \cr
    & =    4 \bigl\| P_{n_{k+1}} (I - Q_A) (x) \bigr\|_J ^2
      \leq 16 \|x\|_J^2, \cr
}$$
so $S(x) = ( \< x^*_k, T_k(x) \> )_{k\in\N}$ defines a bounded operator
$\map S\JJphi{\ell_2}$ which maps, for every $j\in\N$,
$$S(z_j) = ( \< x^*_k, T_k(z_j) \> )_{k\in\N} =
  ( \< x^*_k, \delta_{jk} z_j \> )_{k\in\N} = e_j$$
As the operator $\ell_2\maparrow\JJphi$ that takes $(\alpha_n)_{n\in\N}$ to
$\sum_{n=1}^\infty \alpha_n z_n$ (hence each $e_j$ to~$z_j$) is bounded
by Proposition~\ref{J-upper}, it follows that $Z = [z_k : k\in\N]$ is
isomorphic to~$\ell_2$ and complemented in~$\JJphi$.

Finally, in the spirit of the principle of small perturbations
\cite{bessaga-pelcz},\break let $\map K\JJphi\JJphi$ be the operator defined
as $K(x) = $\break$\sum_{k=1}^\infty \< x^*_k, T_k(x) \> ({x_k - z_k})$; then
$\sum_{k=1}^\infty \|x^*_k\| \|T_k\| \|x_k - z_k\| < 1$,
so $K$ is well defined and $U = I + K$ is an isomorphism on~$\JJphi$, where
$K(z_k) = 
x_k - z_k$ and $U(z_k) = x_k$ for every $k\in\N$, and then $U(Z) =
[x_k : k\in\N] \subseteq M$ is a copy of~$\ell_2$ complemented in~$\JJphi$.
\end{Proof}

\begin{Cor}
\label{J-subprojective}
$\Jphi$ is subprojective if and only if each $X_n$ is subprojective.
\end{Cor}

The fact that any closed infinite-dimensional subspace of~$\Jphi$ contains
another infinite-dimensional subspace complemented in~$\JJphi$, and not
merely in~$\Jphi$, allows to show that $\JJphi$ is subprojective when
additional conditions on $\JJphi/\Jphi$ are met.

\begin{Lemma}
\label{sub-3sp-lemma}
Let $X$ be a Banach space and let $M$ and $N$ be closed 
subspaces of~$X$ such that $M \cap N = 0$ and $M + N$ is not closed.
Then there exists an automorphism $\map UXX$ such that $U(M) \cap N$ is
infinite-dimensional.
\end{Lemma}

\begin{Proof}
Take normalised sequences $(x_n)_{n\in\N}$ in~$M$ and $(y_n)_{n\in\N}$ in~$N$
such that $\|x_n - y_n\| < 2^{-n}$ for every $n\in\N$. Since any weak
cluster point of $(x_n)_{n\in\N}$ must be in $M\cap N = 0$, by passing
to a subsequence \cite[Theorem~1.5.6]{kalton-albiac} we can assume that
$(x_n)_{n\in\N}$ is a basic sequence and that there exists a sequence
$(x^*_n)_{n\in\N}$ in~$X^*$ such that $\<x^*_i,x_j\> = \delta_{ij}$ for
every $i$, $j\in\N$ and $\sum_{n=1}^\infty \|x^*_n\| \, \|x_n - y_n\| < 1$.
Then $K(x) = \sum_{n=1}^\infty \<x^*_n,x\> (x_n - y_n)$ defines an operator
$\map KXX$ with $\|K\| < 1$ and $U = I - K$ is an automorphism on~$X$
that maps $U(x_n) = y_n$ for every $n\in\N$, so $U(M) \cap N$ is
infinite-dimensional.
\end{Proof}

\begin{Prop}
\label{JJ-subprojective}
If $\Jphi$ and $\JJphi/\Jphi$ are both subprojective,
then $\JJphi$ is subprojective.
\end{Prop}

Note that subprojectivity is not a three-space property, in general
\cite[Proposition 2.8]{oikhberg-spinu}.

\begin{Proof}
Let $M$ be a closed infinite-dimensional subspace of~$\JJphi$.
If $M\cap\Jphi$ is infinite-dimensional, then it contains another
infinite-dimensional subspace complemented in~$\JJphi$ by
Theorem~\ref{J-main}.

Otherwise, if $M\cap\Jphi$ is finite-dimensional, we can assume that
$M\cap\Jphi = 0$ by passing to a further subspace if necessary.
If $M + \Jphi$ is closed, let $\map{Q}{\JJphi}{\JJphi/\Jphi}$ be the
natural quotient operator induced by~$\Jphi$;
then $Q|_M$ is an embedding into~$\JJphi/\Jphi$, which is subprojective,
so $M$ contains an infinite-dimensional subspace complemented in~$\JJphi$
by Proposition~\ref{subproj-helper}.

We are left with the case where $M\cap\Jphi = 0$ and $M + \Jphi$ is not
closed. By Lemma~\ref{sub-3sp-lemma}, there exists an automorphism
$\map U\JJphi\JJphi$ such that $U(M) \cap \Jphi$ is infinite-dimensional.
Let $N$ be an infinite-dimensional subspace of $U(M) \cap \Jphi$ complemented
in~$\JJphi$, which exists again by Theorem~\ref{J-main}.
then $U^{-1}(N) \subseteq M$ and is still complemented in~$\JJphi$.
\end{Proof}

\begin{Prop}
\label{JJ-subprojective-SS}
If $\Jphi$ is subprojective and the quotient map\break $\JJphi \maparrow
\JJphi/\Jphi$ is strictly singular, then $\JJphi$ is subprojective.
\end{Prop}

\begin{Proof}
Let $M$ be a closed infinite-dimensional subspace of~$\JJphi$. If $M\cap\Jphi$
is infinite-dimensional, then it contains another infinite-dimensional
subspace complemented in~$\JJphi$ by Theorem~\ref{J-main}.

Otherwise, if $M\cap\Jphi$ is finite-dimensional, we can assume that
$M\cap\Jphi = 0$ by passing to a further subspace if necessary, and then
the strict singularity of the quotient map implies that $M + \Jphi$
cannot be closed. By Lemma~\ref{sub-3sp-lemma}, there exists an automorphism
$\map U\JJphi\JJphi$ such that $U(M) \cap \Jphi$ is infinite-dimensional.
Let $N$ be an infinite-dimensional subspace of $U(M) \cap \Jphi$ complemented
in~$\JJphi$, which exists again by Theorem~\ref{J-main};
then $U^{-1}(N) \subseteq M$ and is still complemented in~$\JJphi$.
\end{Proof}

A Banach space~$X$ is called \emph{superprojective} if every closed
infinite-codimensional subspace of~$X$ is contained in an
infinite-codimensional subspace complemented in~$X$~\cite{whitley}
(see also~\cite{superprojective}).
There is a perfect duality between subprojectivity and superprojectivity
for reflexive Banach spaces, in that a reflexive Banach space is
subprojective if and only if its dual is superprojective. This does not
extend to the non-reflexive case: $c_0$ is subprojective while $\ell_1$
is not superprojective \cite{whitley}, and $\ell_1$ is subprojective
but it has a hereditarily indecomposable predual \cite{argyros-haydon}.
It is still unclear whether the superprojectivity of a Banach space
implies the subprojectivity of its dual or predual.
In light of Corollary~\ref{J-subprojective} and
Proposition~\ref{JJ-subprojective}, it makes sense
to ask whether $\Jphi$ would be superprojective if so were every $X_n$, or
what additional conditions on the $X_n$ would be required.

\begin{Question}
If every $X_n$ is superprojective, is $\Jphi$ superprojective?
\end{Question}

The proof of Theorem~\ref{J-main} relies on the fact that sequences
in~$\JJphi$ with skipped disjoint supports generate (complemented) copies
of~$\ell_2$, due to Propositions \ref{J-upper} and~\ref{J-lower}. If the
supports are disjoint but not skipped then only Proposition~\ref{J-upper}
applies, so other subspaces are possible. However, the $J$-sum never
introduces copies of~$\ell_p$ for $p > 2$. This is already proved in
\cite[Lemma~5]{ostrovskii} when every $X_n$ is finite-dimensional,
but this condition can be relaxed to the case where no $X_n$ contains
a copy of~$\ell_p$. If the construction of~$\Jphi$ is done with an
arbitrary $1 < q < \infty$ instead of~$2$, it is also known that $\JJphi$
does not contain a copy of any $\ell_p$ for $p > q$ if every $X_n$ is
finite-dimensional \cite[Lemma~3]{quojections}; as mentioned in
Remark~\ref{not-2}, the results below extend to this case as well,
so no space $\JJphi$ constructed using any $1 < q < \infty$ contains
a copy of~$\ell_p$ for $p > q$ unless one of the $X_n$ already does.

\begin{Lemma}
\label{J-q-p-lemma}
Let $(m_k)_{k\in\N}$ be a strictly increasing sequence of positive integers,
let $m_0 = 0$ and let $x^k \in R((I-P_{m_{k-1}})Q_{m_k})$ for every $k\in\N$.
Then $2 \bigl\|\sum_{i=1}^k x^i\bigr\|_J^2 \geq \sum_{i=1}^k \|x^i\|_J^2$.
\end{Lemma}

\begin{Proof}
For every $k\in\N$, write $x^k = (x^k_n)_{n\in\N_0}$ and note that
$P_{m_{k-1}}Q_{m_k} = Q_{m_k}P_{m_{k-1}}$, as $m_{k-1} < m_k$, so $x^k \in
R((I-P_{m_{k-1}})Q_{m_k}) = R(Q_{m_k}(I-P_{m_{k-1}}))$. As a consequence,
by Proposition~\ref{J-interval} and Lemma~\ref{J-Q-S},
$$\sup_{S\in\psf} \rho(x^k,S) = \sup_{S\in\psf}
  \rho \bigl (x^k, \{m_{k-1}\} \cup (S \cap (m_{k-1},m_k]) \bigr),$$
so there exists some $S\in\psf$ such that $m_{k-1}\in S \subseteq
[m_{k-1},m_k]$ and $\rho(x^k,S)^2 = 2\|x^k\|_J^2$; since $x^k_{m_{k-1}} = 0$,
letting $p = \max S \leq m_k$, we have
$$2\|x^k\|_J^2 = \rho(x^k,S)^2 = \sigma(x^k,S)^2 + \|x^k_p\|_p^2 =
  \sigma(x^k,S)^2 + \sigma(x^k,\{m_{k-1},p\})^2,$$
so there exists some $S_k\in\psf$ (either $S$ or $\{m_{k-1},p\}$)
such that $m_{k-1}\in S_k \subseteq [m_{k-1},m_k]$ and
$\sigma(x^k,S_k)^2 \geq \|x^k\|_J^2$.

Fix now $k\in\N$ and consider $\sigma \bigl( \sum_{i=1}^k x^i, S_j \bigr) ^2$
for a given $j\in\{1,\ldots,k\}$. For any two elements $n < p \in S_j$,
$$\phi_n^p \biggl( \sum_{i=1}^k x^i_n \biggr) - \sum_{i=1}^k x^i_p =
  \sum_{i=1}^k \bigl( \phi_n^p(x^i_n) - x^i_p \bigr),$$
where, for $i < j$, we have $x^i \in R(Q_{m_i})$, so
$m_i \leq m_{j-1} \leq n < p$ and
$$\phi_n^p(x^i_n) - x^i_p =
  \phi_n^p(\phi_{m_i}^n(x^i_{m_i})) - \phi_{m_i}^p(x^i_{m_i}) = 0
\vadjust{\break}$$
and, for $i > j$, we have $x^i \in R(I-P_{m_{i-1}})$, so
$n < p \leq m_j \leq m_{i-1}$ and
$$\phi_n^p(x^i_n) - x^i_p = \phi_n^p(0) - 0 = 0$$
so only $i = j$ matters and
$$\phi_n^p \biggl( \sum_{i=1}^k x^i_n \biggr) - \sum_{i=1}^k x^i_p =
  \sum_{i=1}^k \bigl( \phi_n^p(x^i_n) - x^i_p \bigr) =
  \phi_n^p(x^j_n) - x^j_p,$$
hence
$$\eqalign{
  \rho   \biggl( \sum_{i=1}^k x^i, \bigcup_{j=1}^k S_j \biggr) ^2
    & \geq  \sigma \biggl( \sum_{i=1}^k x^i, \bigcup_{j=1}^k S_j \biggr) ^2
      \geq  \sum_{j=1}^k \sigma \biggl( \sum_{i=1}^k x^i, S_j \biggr) ^2 \cr
    & =     \sum_{j=1}^k \sigma(x^j,S_j)^2
      \geq  \sum_{j=1}^k \|x^j\|_J^2. \cr
}$$
\end{Proof}

\begin{Prop}
\label{J-q-p}
Let $p > 2$.
If no $X_n$ contains a copy of~$\ell_p$, then neither does $\JJphi$.
\end{Prop}

\begin{Proof}
Assume otherwise and let $(x_k)_{k\in\N}$ be a sequence in~$\JJphi$
equivalent to the unit vector basis of~$\ell_p$.

Then, for every $j$, $m\in\N_0$ and $\eps > 0$, $R(P_m) \simeq
\oplus_{i=1}^m X_i$ does not contain copies of~$\ell_p$, as containing copies
of~$\ell_p$ is a three-space property~\cite{3space}, so $P_m|_{[x_i:j<i]}$ is
strictly singular. Therefore, there exist $k\in\N$ with $k > j$ and
$y \in [x_i:j<i\leq k]$ such that $\|y\|_J = 1$ and $\|P_m(y)\|_J < \eps/2$,
and then there exists $n > m$ such that $\|(I-Q_n)(I-P_m)(y)\|_J < \eps/2$
by Proposition~\ref{J-Qlim}, which leads to $\|y - Q_n(I-P_m)(y)\|_J \leq
\|P_m(y)\|_J + \|y - P_m(y) - Q_n(I-P_m)(y)\|_J < \eps$.

By induction, given a sequence $(\eps_n)_{n\in\N}$ of positive real numbers
and starting with $m_0 = 0$, there exist a strictly increasing sequence
$(m_n)_{n\in\N}$ of positive integers and a sequence $(y_n)_{n\in\N}$ of
norm-one elements of~$\JJphi$ that is a blocking of $(x_n)_{n\in\N}$ such
that $\|y_n - Q_{m_n}(I-P_{m_{n-1}})(y_n)\|_J < \eps_n$ for every $n\in\N$;
in particular, $(y_n)_{n\in\N}$ is still equivalent to the unit vector basis
of~$\ell_p$.
Choosing the right sequence $(\eps_n)_{n\in\N}$, we can then assume that
$(z_n)_{n\in\N} = \bigl( Q_{m_n}(I-P_{m_{n-1}})(y_n) \bigr)_{n\in\N}$
is also still equivalent to the unit vector basis of~$\ell_p$, where each
$z_n \in R \bigl( Q_{m_n} (I-P_{m_{n-1}}) \bigr) =
R \bigl( (I-P_{m_{n-1}}) Q_{m_n} \bigr)$. But then, for any $n\in\N$,
we would have
$$2 \biggl\|\sum_{i=1}^n z_i\biggr\|_J^2 \geq \sum_{i=1}^n \|z_i\|_J^2$$
by Lemma~\ref{J-q-p-lemma}, which is not possible for a sequence equivalent
to the unit vector basis of $\ell_p$ with $p > 2$.
\end{Proof}

\break

For example, let $(p_n)_{n\in\N}$ be a strictly increasing sequence
in $(2,\infty)$ and let $p = \lim_n p_n$. Taking $X_n = \ell_{p_n}$ and
letting $\map{\phi_n}{\ell_{p_n}}{\ell_{p_{n+1}}}$ be the natural inclusion,
then $\JJphi/\Jphi$ is isomorphic to~$\ell_p$ (or $c_0$ if $p = \infty$),
as $\smash{\|x\|_{p_n} \tendsto_n \|x\|_p}$, so the quotient map $\JJphi
\maparrow \JJphi/\Jphi$ is strictly singular by Proposition~\ref{J-q-p}
and $\JJphi$ is subprojective by Proposition~\ref{JJ-subprojective-SS}.

\end{section}

\begin{section}{The dense subspace chain case}

\def\JX{{J(X_n)}}
\def\JJX{{\JJ(X_n)}}

Let $(Z, \|\cdot\|)$ be a Banach space. In this section, we will assume that
each $X_n$ is a closed subspace of~$Z$ such that $\bigcup_{n\in\N} X_n$
is dense in~$Z$ and $X_n \subseteq X_{n+1}$ for every $n\in\N_0$,
each $\|\cdot\|_n$ is the restriction of~$\|\cdot\|$ to~$X_n$ and each
$\phi_n$ is the natural inclusion of $X_n$ into~$X_{n+1}$.
We will write $\JX = \Jphi$ and $\JJX = \JJphi$ for this particular case.

Under these conditions, the mapping $\map{S}{\JJX}{Z}$ defined as
$S \bigl( (x_n)_{n\in\N_0} \bigr) = \lim_n x_n$ is a surjective bounded
linear operator with kernel~$\JX$ and the induced quotient map
$\JJX/\JX \maparrow Z$ is an isometry \cite[Corollary 1.3]{Jsum}.

\begin{Prop}
If $Z$ is subprojective, then $\JJX$ is subprojective.
\end{Prop}

\begin{Proof}
Each $X_n \subseteq Z$ is subprojective \cite[Lemma~3.1]{whitley}
and $\JJX/\JX$ is isomorphic to~$Z$, hence subprojective, so $\JJX$
is subprojective by Proposition~\ref{JJ-subprojective}.
\end{Proof}

Given that the construction of~$\JX$ is made so that the quotient $\JJX/\JX$
is isometric to~$Z$, it is natural to ask whether the quotient is
complemented, that is, whether $\JJX = \JX\oplus(\JJX/\JX) = \JX\oplus Z$.
For some spaces~$Z$, this will be impossible; for instance, $Z = c_0$ as a
subspace of $\JJX \equiv \JX{}^{**}$ would make $\ell_\infty$ embed
into~$\JJX$, which is separable ($\Omega(X_n)$ is dense in~$\JJX$).
On the other hand, for $Z = \ell_1$ the complementation is immediate,
as $\ell_1$ is always complemented as a quotient.
For other spaces, we have the following.

\begin{Prop}
Let $(Z_n)_{n\in\N}$ be an unconditional Schauder decomposition of~$Z$
that satisifies a lower $2$-estimate and such that $X_n = \oplus_{i=1}^n Z_i$
and let $(T_n)_{n\in\N}$ be the sequence of aggregate projections associated
with $(Z_n)_{n\in\N}$ (so $R(T_n) = X_n$ and $\smash{T_n(z) \tendsto_n z}$).
Then the linear map $\map{T}{Z}{\JJX}$ defined as $T(z) = (T_n(z))_{n\in\N}$
is well defined and an embedding into~$\JJX$,
and $\JJX = \JX \oplus R(T) \simeq \JX \oplus Z$.
\end{Prop}

\begin{Proof}
Let $K < \infty$ be the suppresion constant of $(Z_n)_{n\in\N}$ and
let $M < \infty$ be such that $\sum_{i=1}^k \|z_i\|^2 \leq
M \bigl\| \sum_{i=1}^k z_i \bigr\|^2$ if $(z_i)_{i=1}^k$ are disjoint
with respect to $(Z_n)_{n\in\N}$. Let $z\in Z$ and
$S = \{\, p_0 < \cdots < p_k \,\} \in \psf$; then
$$\eqalign{
  \rho(T(z),S)^2 &= \sum_{i=1}^k \bigl\| T_{p_{i-1}}(z)
        - T_{p_i}(z) \bigr\| ^2 + \|T_{p_k}(z)\|^2 \cr
    &\leq M \biggl\| \sum_{i=1}^k \bigl( T_{p_{i-1}}(z)
        - T_{p_i}(z) \bigr) \biggr\| ^2 + \|T_{p_k}(z)\|^2 \cr
    &=    M \bigl\| T_{p_0}(z) - T_{p_k}(z) \bigr\| ^2 + \|T_{p_k}(z)\|^2 \cr
    &\leq M K^2 \|z\|^2 + K^2 \|z\|^2 = (M{+}1) K^2 \|z\|^2 \cr
}$$
so $T(z) \in \JJX$ and $\|T(z)\|_J \leq \sqrt{(M{+}1)/2} K \|z\|$.
If we define the operator $\map{S}{\JJX}{Z}$ given by
$S \bigl( (x_n)_{n\in\N} \bigr) = \lim_n x_n$ as before,
then $ST$ is clearly the identity on~$Z$, so $TS$ is a projection on~$\JJX$
with kernel~$\JX$ and range $R(T)$.
\end{Proof}

This applies, for instance, to $Z = \ell_p$, for $1\leq p\leq 2$, with its
canonical decomposition, which is unconditional and satisfies a lower
$2$-estimate \cite[Theorem~1.f.7]{classical2}, or to $Z = L_p$ for
$1\leq p\leq 2$ with the decomposition associated to the Haar basis
(or any other unconditional basis) \cite[p.~128]{alspach-odell}.
On the other hand, for $Z = \ell_p$ with $p > 2$ and its canonical
decomposition (or any other space that contains a copy of~$\ell_p$ for
$p > 2$ and a decomposition where no $X_n$ does), $\JJX$ does not contain
copies of~$\ell_p$ by Proposition~\ref{J-q-p}, so $\JJX$ cannot contain
$\JJX/\JX \equiv Z$ as a subspace.

\end{section}

\end{document}